\documentstyle[12pt]{article}
\author{Youssef ALAOUI}
\title{A generalization of the Levi problem with singularities}
\date{}
\newcommand{\reels}{\mbox{I}\!\!\!\mbox{R}}

\newcommand{\complexes}{\mbox{I}\!\!\!\mbox{C}}

\newtheorem{th}{theorem}
\newtheorem{lm}{lemma}
\begin{document}
\maketitle
\setcounter{page}{1}
\noindent
{\large Summary}. In this article, we prove that if $X$
is a Stein space and $\Omega\subset X$ an increasing sequence
of $q$-complete open subsets, then $\Omega$ is $q$-complete.
\\
$2000$ MS Classification numbers: 32E10, 32E40.\\
\\
\\
{\bf 1. Introduction}\\
\\
\hspace*{.1in}It is known from a classical theorem due to Benke and
 Stein $[1]$ that if
$D_{1}\subset D_{2}\subset\cdots\subset D_{n}\subset\cdots$ is an
 increasing
sequence of Stein open sets in $\complexes^{n}$, then their union
$\displaystyle\bigcup_{j\geq 1}D_{j}$ is Stein.\\
\hspace*{.1in}This result has been extended to open sets in Stein
 manifolds. But
nothing is known about this if $(D_{j})_{j\geq 1}$ is a family of a
Stein complex space $X$, even when $X$ has isolated singularities.\\
\hspace*{.1in}The main new result of this article concerns a
 generalization
of this theorem for families of $q$-complete open sets in Stein spaces.
By using a description of $q$-convex functions in terms of linear sets,
 introduced by
M. Peternell $[4]$, we show that if $X$ is a Stein space and
 $(\Omega)_{j\geq 1}$
an increasing sequence of $q$-complete open subsets of $X$, then
$\Omega=\displaystyle\bigcup_{j\geq 1}\Omega_{j}$ is $q$-complete.

For this, we prove that there exist a linear set ${\mathcal{M}}$
of codimension $\leq q-1$ over $\Omega$ and an exhaustion function
 $\psi\in C^{\infty}(\Omega)$
such that every point of $\Omega$ has a neighborhood $U$ on which there
 exist
finitely many smooth functions $\psi_{1},\cdots, \psi_{l}$ which are
 $1$-convex with respect
to ${\mathcal{M}}|_{U}$, and such that
$$\psi|_{U}=max(\psi_{1},\cdots, \psi_{l})$$
\\
\\
{\bf 2. Preliminaries}\\
\\
\hspace*{.1in}Let $X$ be an analytic complex space. A linear set over
 $X$
is a subset ${\mathcal{M}}$ of $TX=\displaystyle\bigcup_{x\in
 X}T_{x}X$,
where $T_{x}X$ is the Zariski tangent space of $X$ at $x$, such that
 for every
point $x\in X$, ${\mathcal{M}}_{x}={\mathcal{M}}\cap T_{x}X$
is a complex vector subspace of $T_{x}X$.\\
\hspace*{.1in}We recall that a smooth real-valued function $\phi$ on
 $X$ is said to be\\
\hspace*{.1in}a) $q$-convex, if every point $x\in X$ has an open
 neighborhood $U$ isomorphic
to a closed analytic set in a domain $D\subset \complexes^{n}$ such
 that the restriction
$\phi|_{U}$ has an extension $\tilde{\phi}\in C^{\infty}(D)$ whose Levi
 form $L(\tilde{\phi},z)$
has at most $q-1$ negative or zero eingenvalues at any point $z$ of
 $D$.\\
\hspace*{.1in}b) weakly $1$-convex with respect to the linear set
 ${\mathcal{M}}$, if for every
point $x\in X$, there exist a holomorphic embedding $i: U\rightarrow
 \hat{U}$,
where $U\subset X$ is an open neighborhood of $x$ and $\hat{U}$ an open
 subset of some
$\complexes^{N}$, and a extension $\tilde{\phi}\in C^{\infty}(\hat{U},
 \reels)$
of $\phi|_{U}$ such that for every  $\xi\in {\mathcal{M}}_{x}$,
the levi form $L(\tilde{\phi}, i(x))i_{*, x}(\xi)$ is
 positive-semidefinite,
where $i_{*, x}: T_{x}X\rightarrow \complexes^{n}$ denotes the
 differential map of $i$ at $x$.\\
\hspace*{.1in}c) $1$-convex with respect to ${\mathcal{M}}$, if for
 every point $x\in X$,
there exists an open neighborhood $U\subset X$ of $x$ and a $1$-convex
 function
$\theta$ on $U$ such that $\phi-\theta$ is weakly $1$-convex with
 respect to
${\mathcal{M}}|_{U}=\displaystyle\bigcup_{x\in U}{\mathcal{M}}_{x}$.\\
\\
\hspace*{.1in}Let $\Omega$ be an open set in $X$. We denote by
${\mathcal{B}}(\Omega, {\mathcal{M}}|{\Omega})$
the set of all continuous functions on $\Omega$ which are locally the
 supremum
of finitely many $1$-convex functions with respect to
 ${\mathcal{M}}|{\Omega}$
and we define
 $codim_{\Omega}{\mathcal{M}}=sup_{x\in\Omega}{\mathcal{M}}_{x}$.\\
\\
\hspace*{.1in}The space $X$ is said to be $q$-complete, if there exists
a function $\phi\in C^{\infty}(X, \reels)$ which is $q$-convex on $X$
and such that $\phi$ is an exhaustion function i.e.\\ $\{x\in X:
 \phi(x)<c\}$
is relatively compact in $X$ for any $c\in\reels$.\newpage
\noindent
{\bf 3. Proof of theorem $1$.}\\
\\
\hspace*{.1in}We shall prove theorem $1$ using the following result of
M. Peternell $[6]$:
\begin{lm}{Let $X$ be a complex space and $\phi\in C^{\infty}(X,
 \reels)$
a $q$-convex function. Then there exists a linear set ${{\mathcal{M}}}$
 of codimension $\leq q-1$
over $X$ such that $\phi$ is $1$-convex with respect to
 ${{\mathcal{M}}}$.}
\end{lm}
\hspace*{.1in}Let $X$ be a Stein space of dimension $n$ and
 $\Omega\subset X$ an open subset which is
the union of an increasing sequence $\Omega_{1}\subset
 \Omega_{2}\subset\cdots\subset \Omega_{n}\subset\cdots$
of $q$-complete open sets in $X$. Let $\phi_{\nu}:
 \Omega_{\nu}\rightarrow ]0, +\infty[$
be a smooth $q$-convex exhaustion function on $\Omega_{\nu}$,
 ${\mathcal{M}}_{\nu}\subset T\Omega_{_{\nu}}$
a linear set of codimension $\leq q-1$ such that $\phi_{\nu}$ is
 $1$-convex
with respect to ${\mathcal{M}}_{\nu}$ over $\Omega_{\nu}$,
$\nu\geq 1$, and let $(d_{\nu})_{\nu\geq 1}$
be a sequence with $d_{\nu}<d_{\nu+1}$, and $Sup d_{\nu}=+\infty$. One
 may assume that
if $\Omega'_{\nu}=\{x\in \Omega_{\nu}: \phi_{\nu}(x)<d_{\nu}\}$,
then $\Omega'_{\nu}\subset\subset \Omega'_{\nu+1}$.
\begin{lm}{-There exist for each $\nu\geq 1$ an exhaustion function\\
 $\varphi_{\nu}\in C^{\infty}(\Omega_{\nu})$
which is $q$-convex in a neighborhood
of $\overline{\Omega'}_{\nu}\setminus\Omega'_{\nu-1}$, a locally finite
 covering $(U_{\nu})_{\nu\geq 1}$
of $\Omega$ by open sets $U_{\nu}\subset \Omega_{\nu+1}$, a linear set
 ${\mathcal{M}}$
of codimension $\leq q-1$ over $\Omega$ and constants $c_{\nu}\in
 \reels$, $\nu\geq 1$,
with the following properties:\\
(a) For each $\nu\geq 1$ there exists $\psi_{\nu}\in
 {\mathcal{B}}(U_{\nu}, {\mathcal{M}}|_{U_{\nu}})$
with $\psi_{\nu}=\psi_{\nu-1}$ on\\ $\{x\in U_{\nu}:
 \varphi_{\nu+1}(x)<c_{\nu}\}\cap U_{\nu-1}$.\\
(b) For every index $\nu\geq 1$, there is $\varepsilon_{\nu}>0$ such
 that\\
$\Omega'_{\nu-1}\setminus\overline{\Omega'}_{\nu-2}\subset  \{x\in
 U_{\nu}: \varphi_{\nu+1}(x)<c_{\nu}-\varepsilon_{\nu}\}$\\
and $\{x\in U_{\nu}:
 \varphi_{\nu+1}(x)<c_{\nu}+\varepsilon_{\nu}\}\subset U_{\nu-1}$.}
\end{lm}
Proof. There exists a $C^{\infty}$ exhaustion function
 $\varphi_{\nu+1}$ on $\Omega_{\nu+1}$
which is $q$-convex in a neighborhood of
 $\overline{\Omega'}_{\nu+1}\setminus\Omega'_{\nu}$ such that, if\\
$m_{\nu+1}=Min_{\{\phi_{\nu}\geq d_{\nu}\}\cap
 \overline{\Omega'}_{\nu+1}} \varphi_{\nu+1}$
and $M_{\nu+1}=Max_ {\{\phi_{\nu-1}\leq d_{\nu-1}\}}\varphi_{\nu+1}$,
then\\ $m_{\nu+1}>M_{\nu+1}$.
In fact, let $\theta_{\nu}\geq 0$ be a $C^{\infty}$ function with
 compact support
in $\Omega_{\nu+1}\setminus\overline{\Omega'}_{\nu-1}$ such that
 $\theta_{\nu}=1$
on a neighborhood of\\
 $\overline{\Omega'}_{\nu+1}\setminus\Omega'_{\nu}$.
Let $\xi$ be a point of $\partial{\Omega'}_{\nu-1}$ such that
$\phi_{\nu+1}(\xi)=Max_{\{\phi_{\nu-1}=d_{\nu-1}\}}\phi_{\nu+1}$.
Then it is clear that
$$\varphi_{\nu+1}=\phi_{\nu+1}+\phi_{\nu+1}(\xi)\theta_{\nu}$$
satisfies the requirements.\\
\hspace*{.1in}We now put
$$U_{1}=\Omega'_{2}, \ \ and \ \
 U_{\nu}=(\Omega'_{\nu+1}\setminus\overline{\Omega'}_{\nu-2})
 \ \ for \ \  \nu\geq 2.$$
 Since $(U_{\nu})_{\nu\geq 1}$ is a locally finite covering of $\Omega$
and, for every $\nu\geq 1,$ $\phi_{\nu+1}$ is $1$-convex with respect
 to
${{\mathcal{M}}}_{\nu+1}$ over $U_{\nu}$,
it follows from ($[3]$, lemma $2$) that there exists a linear set
${{\mathcal{M}}}$ over $\Omega$ of codimension $\leq q-1$ such that
 each $\phi_{\nu+1}$
is $1$-convex with respect to ${{\mathcal{M}}}$ over $U_{\nu}$.
Moreover, if we set
$$c'_{\nu}=Inf\{\varphi_{\nu+1}(x), x\in
 (\overline{\Omega'}_{\nu+1}\setminus\Omega'_{\nu})\},$$
then
$$(\Omega'_{\nu-1}\setminus\overline{\Omega'}_{\nu-2})\subset \{x\in
 U_{\nu}: \varphi_{\nu+1}(x)<c'_{\nu}\}\subset
(\Omega'_{\nu}\setminus\overline{\Omega'}_{\nu-2})\subset U_{\nu-1}.$$
\hspace*{.1in}Choose $\varepsilon_{\nu}>0$ such that $M_{\nu+1}<c'_{\nu}-2\varepsilon_{\nu},$ and put $c_{\nu}=c'_{\nu}-\varepsilon_{\nu}$.
Then $\Omega'_{\nu-1}\backslash\Omega'_{\nu-2}\subset\{x\in U_{\nu}: \varphi_{\nu+1}(x)<c_{\nu}-\varepsilon_{\nu}\}$ and\\
$\{x\in U_{\nu}: \varphi_{\nu+1}(x)<c_{\nu}+\varepsilon_{\nu}\}\subset \Omega'_{\nu}\backslash\overline{\Omega'}_{\nu-2}\subset U_{\nu-1}.$\\
\hspace*{.1in}We may take $\theta_{\nu}$ so that\\
$\{x\in U_{\nu}: c_{\nu}+\frac{\varepsilon_{\nu}}{2}\leq
 \varphi_{\nu+1}(x)\leq c_{\nu}+\varepsilon_{\nu}\}\subset \{\theta_{\nu}=1\}$.\\
Then there exists for each $\nu$ a function $\psi_{\nu}:
 \Omega'_{\nu+1}\rightarrow ]0,+\infty[$
with\\ $\psi_{\nu}\in {\mathcal{B}}(U_{\nu},
 {\mathcal{M}}|_{U_{\nu}})$, and such that
$\psi_{\nu}=\psi_{\nu-1}$ on\\ $\{x\in U_{\nu}:
 \varphi_{\nu+1}(x)<c_{\nu}+\frac{\varepsilon_{\nu}}{2}\}.$\\
\hspace*{.1in}In fact, if $\nu=1$, then it is obvious that
$\psi_{1}=\phi_{2}$ has the required properties
for $U_{1}=\Omega'_{2}$ and $\Omega_{1}=\emptyset,$ since 
$\varphi_{2}=\phi_{2}\geq c'_{1}=c_{1}+\varepsilon_{1}$ on $U_{1}$.\\
\hspace*{.1in}We now assume that $\nu\geq 2$ and, that $\psi_{1},
 \cdots, \psi_{\nu-1}$ have been constructed. let
 $\chi_{\nu}(t)=a_{\nu}(t-c_{\nu}-\frac{\varepsilon_{\nu}}{2})$
where $a_{\nu}$ is a positive constant, and consider the function
$\psi_{\nu}: U_{\nu}\rightarrow ]0, +\infty[$ defined by
\[\psi_{\nu}=\left\{\begin{array}{cc}
\psi_{\nu-1} \ \ on \ \  \{\varphi_{\nu+1}\leq
 c_{\nu}-\varepsilon_{\nu}\}\\
Max(\psi_{\nu-1}, \chi_{\nu}(\varphi_{\nu+1})) \ \ on \ \
  \{c_{\nu}-\varepsilon_{\nu}\leq \varphi_{\nu+1}\leq c_{\nu}+\varepsilon_{\nu}\}\\
a_{\nu}(\phi_{\nu+1}+\phi_{\nu+1}(\xi)-c_{\nu}-\frac{\varepsilon_{\nu}}{2})
 \ \ on \ \  \{\varphi_{\nu+1}\geq c_{\nu}+\varepsilon_{\nu}\}
\end{array}
\right.\]
Since on $U'_{\nu}=\{x\in U_{\nu}:
 \varphi_{\nu+1}(x)<c_{\nu}+\frac{\varepsilon_{\nu}}{2}\}\subset U_{\nu-1}$ we have\\
$\psi_{\nu-1}>0>\chi_{\nu}(\varphi_{\nu+1})$ and
$\psi_{\nu-1}\in {\mathcal{B}}(U_{\nu-1}, {\mathcal{M}}|_{U_{\nu-1}})$,
then\\
$\psi_{\nu}=\psi_{\nu-1}\in {\mathcal{B}}(U'_{\nu},
 {\mathcal{M}}|_{U'_{\nu}})$ on $U'_{\nu}$.
On the other hand, the subset\\
 $\{c_{\nu}+\frac{\varepsilon_{\nu}}{2}\leq \varphi_{\nu+1}\leq c_{\nu}+\varepsilon_{\nu}\}\subset U_{\nu-1}$
is contained in $\{\theta_{\nu}=1\}$, which implies that
$\psi_{\nu}=Max(\psi_{\nu-1},
 \chi_{\nu}(\phi_{\nu+1}+\phi_{\nu+1}(\xi)))$ on\\
$\{c_{\nu}+\frac{\varepsilon_{\nu}}{2}\leq \varphi_{\nu+1}\leq
 c_{\nu}+\varepsilon_{\nu}\}$.
Then clearly the function $\psi_{\nu}$ is well-defined and satisfies
 the required conditions,
if $a_{\nu}$ is taken so that
 $a_{\nu}\frac{\varepsilon_{\nu}}{2}>Max_{\{\varphi_{\nu+1}=c_{\nu}+\varepsilon_{\nu}\}}\psi_{\nu-1}$.
\begin{th}{-Let $X$ be a Stein space and
 $\Omega_{1}\subset\Omega_{2}\subset\cdots\subset\Omega_{n}\subset\cdots$
an increasing sequence of $q$-complete open sets in $X$. Then
 $\Omega=\displaystyle\bigcup_{\nu\geq 1}\Omega_{\nu}$
is $q$-complete.}
\end{th}
Proof. Let $\varphi'_{\nu}$ be the function defined by
\[\varphi'_{\nu}=\left\{\begin{array}{cc}
\psi_{\nu} \ \ on \ \  U_{\nu}\setminus\{x\in U_{\nu}:
 \varphi_{\nu+1}(x)\leq c_{\nu}\}.\\
\psi_{\mu} \ \ on \ \  \{x\in U_{\mu+1}: \varphi_{\mu+2}(x)<
 c_{\mu+1}+\frac{\varepsilon_{\mu+1}}{2}\} \ \ for \ \ \mu\leq \nu-1
\end{array}
\right.\]
\hspace*{.1in}Then clearly $\varphi'_{\nu}\in{\mathcal{B}}(V_{\nu},
 {\mathcal{M}}|_{V_{\nu}})$, where
$$V_{\nu}=(U_{\nu}\setminus\{x\in U_{\nu}: \varphi_{\nu+1}(x)\leq
 c_{\nu}\})\cup
(\bigcup_{\mu\leq \nu-1} \{x\in U_{\mu+1}: \varphi_{\mu+2}(x)<
 c_{\mu+1}+\frac{\varepsilon_{\mu+1}}{2}\}).$$
Note that $(V_{\nu})_{\nu\geq 1}$ is an increasing sequence of open
 subsets
of $\Omega$ such that $U_{\nu}\subset V_{\nu}$ and
 $\Omega=\displaystyle\bigcup_{\nu\geq 1}V_{\nu}$.
Moreover, we have $\varphi'_{\nu}=\varphi'_{\nu-1}$ on\\ $\{x\in
 U_{\mu+1}: \varphi_{\mu+2}(x)<c_{\mu+1}\}$
for all $\mu\leq \nu-1$.\\
\hspace*{.1in}Let now $K$ be a compact set in $\Omega$ and $\nu\geq 2$
 such that
$K\subset \Omega'_{\nu-1}$. Since $\varphi'_{\nu}=\varphi'_{\nu-1}$ on
$K\cap
 (\overline{\Omega'}_{\mu}\setminus\overline{\Omega'}_{\mu-1})\subset \{x\in U_{\mu+1}: \varphi_{\mu+2}(x)<c_{\mu+1}\}$
for all $\mu\leq \nu-1$, then $\varphi'_{\nu}=\varphi'_{\nu-1}$ on $K.$
This implies that the sequence $(\varphi'_{\nu})_{\nu\geq 1}$ is
 stationary on every compact subset of $\Omega$.\\
\hspace*{.1in}Since $\{x\in U_{\nu}:
 \varphi_{\nu+1}(x)<c_{\nu}+\varepsilon_{\nu}\}\subset U_{\nu-1}$,
then\\ $\Omega'_{\nu+1}\setminus \Omega'_{\nu}\subset \{x\in U_{\nu}:
 \varphi_{\nu+1}(x)\geq c_{\nu}+\varepsilon_{\nu}\}$,
and hence $\varphi'_{\nu}=\chi_{\nu}(\varphi_{\nu+1})>\nu+1$ on
 $\Omega'_{\nu+1}\setminus \Omega'_{\nu}$,
if $a_{\nu}$ is chosen so that in addition
 $a_{\nu}\frac{\varepsilon_{\nu}}{2}>\nu+1.$\\
\hspace*{.1in}Furthermore, there exists for
each $\nu\geq 1$ a function $\psi'_{\nu}\in C^{\infty}(X)$ such that\\
(a) $\psi'_{\nu}$ is strictly plurisubharmonic on
$\Omega'_{\nu+2}\setminus(\Omega'_{\nu}\setminus\Omega'_{\nu-1}).$\\
(b) $\psi'_{\nu}>2^{\nu+2}$ in
 $\Omega'_{\nu+2}\setminus\overline{\Omega'}_{\nu}$
but $\psi_{\nu}<\varepsilon_{\nu}$ in $\overline{\Omega'}_{\nu-1},$
where $\varepsilon_{\nu}$ is a small positive constant to be chosen
later in the proof. To see this, we choose an open neighborhood
$\omega_{\nu-1}\subset\subset X$ of the holomorphically-convex hull
of $\overline{\Omega'}_{\nu-1}$ in $X$ such that $\omega_{\nu-1}$
contains $\Omega'_{\nu+2}.$
Then there exists a strictly psh exhaustion function $h\in
 C^{\infty}(X)$
such that $h<0$ in $\overline{\Omega'}_{\nu-1}$ and $h>0$ in
 $X\setminus\omega_{\nu-1}.$
Let $\theta\in C^{\infty}_{o}(X),$ $0\leq \theta\leq 1,$ with
$supp(\theta)\subset\omega_{\nu-1},$ $\theta=0$ on
 $\overline{\Omega'}_{\nu-1}$
and $\theta=1$ on $\overline{\Omega'}_{\nu+2}\setminus\Omega'_{\nu}.$
Choose $\varepsilon>0$ sufficiently small so that $h+\varepsilon<0$
on $\overline{\Omega'}_{\nu-1}$ and set for $x\in X$
$$\psi'_{\nu}(x)=A(h(x)+\theta(x)(-\displaystyle
 Min_{y\in\overline{\Omega'}_{\nu+2}}h(y))+\varepsilon)+\varepsilon_{\nu},$$
where $A>0$ is a positive constant.
Then, if $A$ is big enough, $\psi'_{\nu}$ obviously satisfies
properties (a) and (b). Moreover, a similar proof shows that
there exists for each $\mu\leq \nu-1$ a function
$\psi'_{\mu}\in C^{\infty}(X),$
which is strictly psh on\newpage
\noindent
$\Omega'_{\nu+2}\setminus (\Omega'_{\mu}\setminus\Omega'_{\mu-1})$ and,
such that $\psi'_{\mu}>Max_{y\in
 \overline{\Omega'}_{\mu+1}\setminus\Omega'_{\mu}}\psi'_{\mu+1}(y)$
in\\ $\overline{\Omega'}_{\mu+1}\setminus\Omega'_{\mu},$
$\psi'_{\mu}>2^{\mu+2}$ in
 $\overline{\Omega'}_{\mu+2}\setminus\Omega'_{\mu}$
and
$\psi'_{\mu}<\varepsilon_{\mu}$ in $\overline{\Omega'}_{\mu-1},$
where $\varepsilon_{\mu}$ is chosen sufficiently
small so that $\varphi'_{\mu}>\varepsilon_{\mu}$
on $\overline{\Omega'}_{\mu-1}$ for every $1\leq \mu\leq \nu.$
We now consider the function $\psi''_{\nu}: \Omega'_{\nu+1}\rightarrow
 \reels$
defined by\\
$\psi''_{\nu}=Max(\varphi'_{\nu},\psi'_{\nu},\psi'_{\nu-1},\cdots,\psi'_{1})$.
 Then obviously
$\psi''_{\nu}\in {\mathcal{B}}(V_{\nu}, {\mathcal{M}}|_{V_{\nu}})$.
 Moreover,
since
$$V_{\nu+1}\setminus V_{\nu}\subset U_{\nu+1}\setminus\{x\in U_{\nu+1}:
 \varphi_{\nu+2}(x)\leq c_{\nu+1}\}
\subset \Omega'_{\nu+2}\setminus\overline{\Omega'}_{\nu},$$
it follows that for every $j\geq \nu+1$,
 $\psi''_{j}\geq\psi'_{\nu}>2^{\nu+2}$
on $V_{\nu+1}\setminus V_{\nu}$.\\
\hspace*{.1in}Let now $K\subset \Omega$ be a compact subset and
 $\nu\geq 2$
such that $K\subset \Omega'_{\nu-1}$.
Since $\varphi'_{\nu}>\varepsilon_{\nu}$ on
 $\overline{\Omega'}_{\nu-1},$
then
 $\psi''_{\nu}=Max(\varphi'_{\nu},\psi'_{\nu-1},\psi'_{\nu-2},\cdots,\psi'_{1})$
on $\Omega'_{\nu-1}$. Moreover, Since $\varphi'_{\nu}=\varphi'_{\nu-1}$
 on $K$, then $\psi''_{\nu}=\psi''_{\nu-1}$
on $K$, which implies that the sequence $(\psi''_{\nu})_{\nu\geq 1}$
is stationary on every compact subset of $\Omega$.\\
\hspace*{.1in}This proves that the limit $\psi''$ of $(\psi''_{\nu})$
 is an exhaustion
function on $\Omega$ such that $\psi''\in{\mathcal{B}}(\Omega,
 {\mathcal{M}})$.
Since $Codim{\mathcal{M}}\leq q-1$, then $\Omega$ is $q$-complete.
This follows from the fact that for every $\eta\in C^{o}(\Omega,
 \reels)$, $\eta>0$,
there exists a $q$-convex function $\psi$ on $\Omega$ such that
$$\psi''\leq \psi<\psi''+\eta$$
according to a result of Coltiou and Vajaitu (See $[2]$ and $[5]$).


\begin{thebibliography}{10}



\bibitem{bib1} H. Behnke, K. Stein, Konvergente Folgen Von
 Regularitatsbereichen
and die Meromorphiekonvexitat, Math Ann. 166, 204 216(1938).\\
\\
\bibitem{bib2}M. Coltoiu, Complete locally pluripolar sets, J. reine
 angew. Math.,
412 (1992), 108-112.\\
\\
\bibitem{bib3}M. Coltoiu, V. Vajaitu, On the $n$-completeness of
 covering
spaces with parameters. Math. Z. 237, 815-831 (2001).\\
\\
\bibitem{bib4}M. Peternell, Algebraische Variet\"aten und
 $q$-vollst\"andige
komplexe R\"aume. Math. Z., 200 (1989), 547-581.\\
\\
\bibitem{bib5}V. Vajaitu, Approximation theorems and homology of
 $q$-Runge
pairs in complex spaces, J. reine angew. Math., 449 (1994), 179-199.\\
\end{thebibliography}
\end{document}